\documentclass[10pt]{article}

\usepackage{amssymb}
\usepackage{graphicx}
\usepackage{amsmath}
\usepackage{array}

\newcommand{\A}{{\cal A}}
\newcommand{\B}{{\cal B}}
\newcommand{\C}{{\cal C}}

\newcommand{\I}{{\cal I}}

\newcommand{\T}{{\cal T}}
\newcommand{\U}{{\cal U}}

\newcommand{\Nh}{\widetilde{N}}

\newcommand{\kk}{\Bbbk}
\newcommand{\RR}{\mathbb{R}}

\newcommand{\NN}{\mathbb{N}}

\newtheorem{theorem}{Theorem}[section]
\newtheorem{defin}[theorem]{Definition}
\newtheorem{proposition}[theorem]{Proposition}
\newtheorem{corollary}[theorem]{Corollary}
\newtheorem{lemma}[theorem]{Lemma}
\numberwithin{equation}{section}

\renewcommand{\dim}{\mathop{\rm dim}\nolimits}

\newcommand{\cl}{\mathop{\rm cl}\nolimits}
\newcommand{\semimat}{\mathop{\rm Semimat}\nolimits}
\newcommand{\matid}{\mathop{\rm MatId}\nolimits}
\newcommand{\ideal}{\mathop{\rm Ideal}\nolimits}
\newcommand{\preim}{\mathop{\rm Preim}\nolimits}
\newcommand{\matpreim}{\mathop{\rm MatPreim}\nolimits}
\newcommand{\coext}{\mathop{\rm Coext}\nolimits}
\newcommand{\matcoext}{\mathop{\rm MatCoext}\nolimits}
\newcommand{\pointedmat}{\mathop{\rm Pointedmat}\nolimits}

\newcommand{\rank}{\mathop{\rm rank}\nolimits}

\begin{document}

\pagestyle{plain}

\title{Semimatroids and their Tutte polynomials}
\label{chapter:semimatroids}
\author{Federico Ardila}
\date{}
\maketitle

\begin{abstract}
We define and study \emph{semimatroids}, a class of objects which
abstracts the dependence properties of an affine hyperplane
arrangement. We show that geometric semilattices are precisely the
posets of flats of semimatroids. We define and investigate the
Tutte polynomial of a semimatroid. We prove that it is the
universal Tutte-Grothendieck invariant for semimatroids, and we
give a combinatorial interpretation for its non-negative
coefficients.
\end{abstract}

\section{Introduction.}

The goal of this paper is to define and study a class of objects
called \emph{semimatroids}. A semimatroid can be thought of as a
matroid-theoretic abstraction of the dependence properties of an
affine hyperplane arrangement. Many properties of hyperplane
arrangements are really facts about their underlying matroidal
structure. Therefore, the study of such properties can be carried
out much more naturally and elegantly in the setting of
semimatroids.

The paper is organized as follows. In Section \ref{sec:semimats}
we define semimatroids, and show how we can think of a hyperplane
arrangement as a semimatroid. The following sections provide
different ways of thinking about semimatroids. Section
\ref{sec:modid} shows how a semimatroid ``extends" to a matroid,
and determines a modular ideal inside it. The semimatroid can be
recovered from the matroid and its modular ideal. Section
\ref{sec:preim} describes the close relationship between
semimatroids and strong maps. Semimatroids are described in terms
of elementary preimages and single-element coextensions. Section
\ref{sec:pointed} gives a bijection between semimatroids and
pointed matroids. Section \ref{sec:semilat} gives a new
characterization of geometric semilattices as posets of flats of
semimatroids, extending the classical correspondence between
geometric lattices and simple matroids.

The final sections are geared towards the study of the Tutte
polynomial of a semimatroid. Section \ref{sec:del-cont} defines
the concepts of duality, deletion and contraction. Section
\ref{sec:Tuttesemi} defines the Tutte polynomial, and shows that
it is the unique Tutte-Grothendieck invariant for the class of
semimatroids. Finally, Section \ref{sec:intext} gives a
combinatorial interpretation for the non-negative coefficients of
the Tutte polynomial.
%



It is worth pointing out that Kawahara discovered semimatroids
independently, and described their Orlik-Solomon algebra in
\cite{Ka03}. Las Vergnas's work on the Tutte polynomial of a
quotient map \cite{La99} is also closely related to our work; we
will say more about this in Section \ref{sec:Tuttesemi}.

\section{Semimatroids.}\label{sec:semimats}

\begin{defin} \label{def:semi}
A \emph{semimatroid} is a triple $(S, \C , r_{\C})$ consisting of
a finite set $S$, a non-empty simplicial complex $\C$ on $S$, and
a function $r_{\C}:\C \rightarrow \NN$, satisfying the following
five conditions.

\noindent{\bf (R1)} If $X \in \C$, then $0 \leq r_{\C}(X) \leq
|X|$.

\noindent{\bf (R2)} If $X, Y \in \C$ and $X \subseteq Y$, then
$r_{\C}(X) \leq r_{\C}(Y)$.

\noindent{\bf (R3)} If $X,Y \in \C$ and $X \cup Y \in \C$, then
$r_{\C}(X)+r_{\C}(Y) \geq r_{\C}(X \cup Y) + r_{\C}(X \cap Y)$.

\noindent{\bf (CR1)} If $X,Y \in \C$ and $r_{\C}(X) = r_{\C}(X
\cap Y)$, then $X \cup Y \in \C$.

\noindent{\bf (CR2)} If $X,Y \in \C$ and $r_{\C}(X)< r_{\C}(Y)$,
then $X \cup y \in \C$ for some $y \in Y - X$.
\end{defin}

We call $S$, $\C$ and $r_{\C}$ the \emph{ground set},
\emph{collection of central sets} and \emph{rank function} of the
semimatroid $(S, \C, r_{\C})$, respectively. Sometimes we will
slightly abuse notation and denote the semimatroid $\C$, when its
ground set and rank function are clear. We will denote subsets of
$S$ by upper case letters, and elements of $S$ by lower case
letters.

We will need the fact that semimatroids satisfy a ``local"
version of (R1) and (R2) and a stronger version of (CR1) and
(CR2), as follows.

\medskip

\noindent{\bf(R2')} If $X \cup x \in \C$ then $r_{\C}(X \cup x) -
r_{\C}(X) = 0$ or $1$.

\noindent{\bf (CR1')} If $X,Y \in \C$ and $r_{\C}(X) = r_{\C}(X
\cap Y)$, then $X \cup Y \in \C$ and $r_{\C}(X \cup Y) =
r_{\C}(Y)$.

\noindent{\bf (CR2')} If $X,Y \in \C$ and $r_{\C}(X)<r_{\C}(Y)$,
then $X \cup y \in \C$ and $r_{\C}(X \cup y) = r_{\C}(X)+1$ for
some $y \in Y - X$.

\medskip

\noindent \emph{Proof of} (R2'). From (R2) we know that $r_{\C}(X
\cup x) \geq r_{\C}(X)$. From (R3) we know that $r_{\C}(X \cup x)
- r_{\C}(X) \leq r_{\C}(x) - r_{\C}(\emptyset)$, and this is $0$
or $1$ by (R1). $\Box$

\medskip

\noindent \emph{Proof of} (CR1'). The hypotheses imply that $X
\cup Y \in \C$. Then (R2) says that $r_{\C}(Y) \leq r_{\C}(X \cup
Y)$, while (R3) says that $r_{\C}(Y) \geq r_{\C}(X \cup Y)$.
$\Box$

\medskip

\noindent \emph{Proof of} (CR2'). By applying (CR2) repeatedly, we
see that we can keep on adding elements $y_1, \ldots, y_k$ of $Y$
to the set $X$, until we reach a set $X \cup y_1 \cup \cdots \cup
y_k \in \C$ such that $r_{\C}(X \cup y_1 \cup \cdots \cup y_k) =
r_{\C}(Y)$. Now we claim that $r_{\C}(X \cup y_i) = r_{\C}(X)+1$
for some $i$.

If that was not the case then, since $r_{\C}(X \cup y_1) =
r_{\C}(X)$, (CR1') applies to $X \cup y_1$ and $X \cup y_2$.
Therefore $X \cup y_1 \cup y_2 \in \C$ and $r_{\C}(X \cup y_1 \cup
y_2) = r_{\C}(X \cup y_2) = r_{\C}(X)$. Then (CR1') applies to $X
\cup y_1 \cup y_2$ and $X \cup y_3$, so $X \cup y_1 \cup y_2 \cup
y_3 \in \C$ and $r_{\C}(X \cup y_1 \cup y_2 \cup y_3) =
r_{\C}(X)$. Continuing in this way, we conclude that $X \cup y_1
\cup \cdots \cup y_k \in \C$ and $r_{\C}(X \cup y_1 \cup \cdots
\cup y_k) = r_{\C}(X)$, a contradiction. $\Box$

\medskip

Let us now explain the connection between semimatroids and
hyperplane arrangements. Given a field $\kk$ and a positive
integer $n$, an \emph{affine hyperplane} in $\kk^n$ is an
$(n-1)$-dimensional affine subspace of $\kk^n$.  A
\emph{hyperplane arrangement} $\A$ in $\kk^n$ is a finite set of
affine hyperplanes in $\kk^n$.

A subset (or \emph{subarrangement}) $\B \subseteq \A$ of
hyperplanes is \emph{central} if the hyperplanes in $\B$ have a
non-empty intersection. The \emph{rank function} $r_{\A}$ is
defined for each central subset $\B$ by the equation $r_{\A}(\B) =
n - \dim \cap \B$.

\begin{proposition} \label{prop:arrissemi}
Let $\A$ be an affine hyperplane arrangement in $\kk^n$. Let
$\C_{\A}$ be the collection of central subarrangements of $\A$,
and let $r_{\A}$ be the rank function of $\A$. Then $(\A,
\C_{\A},r_{\A})$ is a semimatroid.
\end{proposition}

\noindent \emph{Proof.} To each hyperplane $H_i \in \A$ we can
associate a vector $v_i \in \kk^n$ and a constant $c_i \in \kk$,
so that $H_i$ is the set of points $x \in \kk^n$ such that
$v_i\cdot x= c_i$, with the usual inner product on $\kk^n$. It is
easy to see that the rank of a central subset $\{H_{i_1}, \ldots,
H_{i_k}\} \in \C_{\A}$ is equal to the rank of the set $\{v_{i_1},
\ldots, v_{i_k}\}$ in $\kk^n$.

From this point of view, axioms (R1), (R2), (R3) are standard
facts of linear algebra applied to the vector space $\kk^n$. We
now check axioms (CR1) and (CR2).

To check axiom (CR1), assume that $X, Y \in \C$ and $r_{\A}(X) =
r_{\A}(X \cap Y)$. Let $A = \cap X$ be the intersection of the
hyperplanes in $X$, and similarly let $B = \cap Y$. Since $X \cap
Y \subseteq X$ and $r_{\A}(X \cap Y) = r_{\A}(X)$, we must have
$\cap(X \cap Y) = \cap X = A$. Also, $X \cap Y \subseteq Y$
implies $\cap(X \cap Y) \supseteq \cap Y = B$. Therefore $A
\supseteq B$, and every hyperplane in $X \cup Y$ contains $B$. It
follows that $X \cup Y \in \C$.

To check axiom (CR2), assume that $X,Y \in \C$ and
$r_{\A}(X)<r_{\A}(Y)$. Let $L_X = \{L_i \, | \, H_i \in X\}$ and
define similarly $L_Y$. Since $\rank(L_Y) > \rank(L_X)$, there
exists a vector $L \in L_Y$, corresponding to a hyperplane $y \in
Y$, which is not in the span of $L_X$. Thus $y$ has a non-empty
intersection with $\cap X$. $\Box$

\medskip

Semimatroids, like matroids, have several equivalent definitions.
In their context, it is possible to talk about flats, independent
sets, spanning sets, bases, circuits, and most other basic matroid
concepts. We will say more about this in Section \ref{sec:intext}.
Until then, we will use the rank function approach of Definition
\ref{def:semi} throughout most of our treatment. We will also need
some facts about the closure approach, which we now present.

\begin{defin}
For a semimatroid $\C = (S, \C, r_{\C})$ and a set $X \in \C$, the
\emph{closure of $X$ in $\C$} is $\cl_{\C}(X) = \{x \in S \, | \,
X \cup x \in \C, r_{\C}(X \cup x) = r_{\C}(X) \}$.
\end{defin}

We will sometimes drop the subscript and write $\cl(X)$ instead
of $\cl_{\C}(X)$ when it causes no confusion.

\begin{proposition}
The closure operator of a semimatroid satisfies the following
properties, for all $X, Y \in \C$ and $x,y \in S$.

\noindent {\bf (CLR1)} $\cl(X) \in \C$ and $r_{\C}(\cl(X)) =
r_{\C}(X)$.

\noindent {\bf (CL1)} $X \subseteq \cl(X)$.

\noindent {\bf (CL2)} If $X \subseteq Y$ then $\cl(X) \subseteq
\cl(Y)$.

\noindent {\bf (CL3)} $\cl(\cl(X)) = \cl(X).$

\noindent {\bf (CL4)} If $X \cup x \in \C$ and $y \in \cl(X \cup
x) - \cl(X)$, then $X \cup y \in \C$ and $x \in \cl(X \cup y)$.

\end{proposition}

\noindent \emph{Proof.} To check (CLR1), let $\cl(X) = \{x_1,
\ldots, x_k\}$. We repeat the argument of the proof of (CL2').
Since $r_{\C}(X \cup x_1) = r_{\C}(X)$, (CR1') applies to $X \cup
x_1$ and $X \cup x_2$, so $X \cup x_1 \cup x_2 \in \C$ and
$r_{\C}(X \cup x_1 \cup x_2) = r_{\C}(X)$. (CR1') then applies to
$X \cup x_1 \cup x_2$ and $X \cup x_3$, so $X \cup x_1 \cup x_2
\cup x_3 \in \C$ and $r_{\C}(X \cup x_1 \cup x_2 \cup x_3) =
r_{\C}(X)$. Continuing in this way, we conclude that $X \cup x_1
\cup \cdots \cup x_k \in \C$ and $r_{\C}(X \cup x_1 \cup \cdots
\cup x_k) = r_{\C}(X)$.

(CL1) is trivial.

To check (CL2), let $x \in \cl(X)$. Then $X \cup x \in \C$ and
$r_{\C}(X \cup x) = r_{\C}(X)$. Applying (CR1') to $X \cup x$ and
$Y$, we conclude that $Y \cup x \in \C$ and $r_{\C}(Y \cup x) =
r_{\C}(Y)$. Therefore $x \in \cl(Y)$.

We know that $\cl(X) \subseteq \cl(\cl(X))$; so to prove (CL3) it
suffices to check the reverse inclusion. Let $x \in \cl(\cl(X))$.
Then $\cl(X) \cup x \in \C$ and $r_{\C}(\cl(X) \cup x) =
r_{\C}(\cl(X)) = r_{\C}(X).$ Therefore, since $\cl(X) \cup x
\supseteq X \cup x \supseteq X$, we have $X \cup x \in \C$ and
$r_{\C}(X \cup x) = r_{\C}(X)$ also; \emph{i.e.}, $x \in \cl(X)$.

Finally, we check (CL4). The assumption that $y \in \cl(X \cup x)$
implies that $X \cup x \cup y \in \C$ and $r_{\C}(X \cup x \cup
y) = r_{\C}(X \cup x) \leq r_{\C}(X)+1$. Since $X \cup y \in \C$,
the assumption that $y \notin \cl(X)$ implies that $r_{\C}(X)+1 =
r_{\C}(X \cup y)$. These two results together give $r_{\C}(X \cup
x \cup y) = r_{\C}(X \cup y)$; \emph{i.e.}, $x \in \cl(X \cup y)$.
$\Box$

\medskip

We will later need the following definitions.

\begin{defin}
A \emph{flat} of a semimatroid $\C$ is a set $A \in \C$ such that
$\cl(A) = A$. The \emph{poset of flats} $K(\C)$ of a semimatroid
$\C$ is the set of flats of $\C$, ordered by containment.
\end{defin}

\section{Modular ideals.} \label{sec:modid}

In the following sections, we will present bijections between the
class of semimatroids and other important classes of objects.
Figure \ref{fig:bijections} at the end of Section \ref{sec:preim}
should be useful in understanding these bijections, and is worth
keeping in mind especially in Sections \ref{sec:modid},
\ref{sec:preim} and \ref{sec:pointed}.

 In this section, we show that a semimatroid
is essentially equivalent to a pair $(M,\I)$ of a matroid $M$ and
one of its modular ideals $\I$.


We start by showing how we can naturally construct a matroid
$M_{\C}$ from a given semimatroid $(S,\C,r_{\C})$, by extending
the rank function $r_{\C}$ from $\C$ to $2^S$.

\begin{proposition}\label{prop:semi-mat}
Let $\C = (S, \C, r_{\C})$ be a semimatroid. For each subset $X
\subseteq S$, let $r(X) = \max\{r_{\C}(Y) \, | \, Y \subseteq X \,
, \, Y \in \C\}$. Then $r$ is the rank function of a matroid
$M_{\C} = (S,r)$.
\end{proposition}

\noindent \emph{Proof.} It is clear, but worth remarking
explicitly, that $r(X) = r_{\C}(X)$ if $X \in \C$. It will be most
convenient to check the three \emph{local} axioms (R1')-(R3') for
the rank function of a matroid \cite{Br86}. Let $X \subseteq S$
and $a,b \in S$ be arbitrary.

\medskip

\noindent {\bf (R1')} $r(\emptyset)=0$.

This is trivial.

\medskip

\noindent {\bf (R2')} $r(X \cup a) - r(X) = 0$ or $1$.

This is easy. It is immediate from the definition that $r(X \cup
a) \geq r(X)$. Now let $r(X \cup a) = r_{\C}(Y)$ for $Y \subseteq
X \cup a$, $Y \in \C$. Then $Y-a \subseteq X$ is also in $\C$, so
$r(X) \geq r_{\C}(Y-a) \geq r_{\C}(Y)-1 = r(X \cup a)-1$.

\medskip

\noindent {\bf(R3')} If $r(X \cup a) = r(X \cup b) = r(X)$, then
$r(X \cup a \cup b) = r(X)$.

This takes more work. Assume that $r(X \cup a) = r(X \cup b) =
r(X) = s$ but $r(X \cup a \cup b) = s+1$. Let $W \subseteq X \cup
a \cup b$, $W \in \C$ be such that $r_{\C}(W) = s+1$. Notice that
$W$ must contain $a$; otherwise we would have $W \subseteq X \cup
b$ and $r_{\C}(W) > r(X \cup b)$. Similarly, $W$ contains $b$. So
let $W = Z \cup a \cup b$; clearly $Z \subseteq X$.

We have $s+1 = r_{\C}(Z \cup a \cup b) \leq r_{\C}(Z \cup a) + 1
\leq r(X \cup a) + 1 = s+1$. Therefore $r_{\C} (Z \cup a) = s$.
Similarly, $r_{\C} (Z \cup b) = s$. Then, by the submodularity of
$r_{\C}$, $r_{\C} (Z) = s-1$.

Now, since $r(X) = s$, we can find $V \subseteq X$, $V \in \C$
such that $r_{\C}(V) = s$. So we have $V , Z \in \C$ with $s =
r_{\C}(V) > r_{\C}(Z) = s-1$. By (CR2'), we can add an element of
$V$ to $Z$ and obtain a set $Y \in \C$ with $X \supseteq Y
\supseteq Z$ such that $r_{\C}(Y)=s$. Notice that $Z \cup a
\subseteq Y \cup a \subseteq X \cup a$ and $r(Z \cup a) = r(X
\cup a) = s$. Thus $r(Y \cup a) = s.$ Similarly, $r(Y \cup b) =
s$ and $r(Y \cup a \cup b) = s+1.$

Now we have $Y , Z \cup a \cup b \in \C$ with $s+1 = r_{\C}(Z \cup
a \cup b) > r_{\C}(Y) = s$. Once again, (CR2') guarantees that we
can add an element of $Z \cup a \cup b$ to $Y$ to obtain an
element of rank $s+1$ in $\C$. But $Z \subseteq Y$, so the only
elements of $Z \cup a \cup b$ which may not be in $Y$ are $a$ and
$b$. Also, we saw that $r(Y \cup a) = r(Y \cup b) = s$. This is a
contradiction. $\Box$

\medskip

The following definitions will be important to us.

\begin{defin}
A pair $\{X,Y\}$ of subsets of $S$ is a \emph{modular pair} of the
matroid $M = (S,r)$ if $r(X) + r(Y) = r(X \cup Y) + r(X \cap Y)$.
\end{defin}

\begin{defin} \cite{Do76}
A \emph{modular ideal} $\I$ of a matroid $M=(S,r)$ is a non-empty
collection of subsets of $S$ satisfying the following three
conditions.

\noindent{\bf (MI1)} $\I$ is a simplicial complex.

\noindent{\bf (MI2)} $\{a\} \in \I$ for every non-loop $a$ of $M$.

\noindent{\bf (MI3)} If $\{X,Y\}$ is a modular pair in $M$ and
$X,Y \in \I$, then $X \cup Y \in \I$.

\end{defin}

\begin{proposition}\label{prop:modideal}
For any semimatroid $(S,\C,r_{\C})$, the collection $\C$ is a
modular ideal of $M_{\C}$.
\end{proposition}

\noindent \emph{Proof.} We denote the rank function of $M_{\C}$
by $r$ and the rank function of $\C$ by $r_{\C}$. Of course,
$r_{\C}$ is just the restriction of $r$ to $\C$.

Axioms (MI1) and (MI2) of a modular ideal are satisfied
trivially. We reformulate (MI3) as follows:

\noindent {\bf (MI3)} \,\, If $A,B,C \subseteq S$ are pairwise
disjoint, $A \cup B, A \cup C \in \C$ and $r(A \cup B \cup C) -
r(A \cup B) = r(A \cup C) - r(A)$, then $A \cup B \cup C \in \C$.

We can assume that $B$ and $C$ are non-empty; if one of them is
empty, the claim is trivial. We prove (MI3) by induction on $|B| +
|C|$.

\medskip

The first case is $|B| + |C| = 2$; let $B = \{b\}$ and $C =
\{c\}$. First assume that $r(A \cup b)$ and $r(A \cup c)$ are
different; say, $r_{\C}(A \cup b) < r_{\C}(A \cup c)$. By (CR2),
we can add some element of $A \cup c$ to $A \cup b$ and obtain a
set in $\C$. This element can only be $c$, so $A \cup b \cup c
\in \C$.

Assume then that $r_{\C}(A \cup b) = r_{\C}(A \cup c) = s$. If
$r_{\C}(A) = s$, (CR1) implies that $A \cup b \cup c \in \C$.
Assume then that $r_{\C}(A) = s-1$, and therefore $r(A \cup b
\cup c) = s+1$ by hypothesis. There is a subset of $A \cup b \cup
c$ in $\C$ of rank $s+1$; since it cannot be contained in $A \cup
b$ or $A \cup c$, it must be of the form $B \cup b \cup c$ for
some $B \subseteq A$. But then we have $r_{\C}(A \cup b) <
r_{\C}(B \cup b \cup c)$. By (CR2), we can add some element of $B
\cup b \cup c$ to $A \cup b$ and obtain a set in $\C$. This
element can only be $c$, so $A \cup b \cup c \in \C$.

\medskip

Having established the base case of the induction, we proceed
with the inductive step. Assume that $|B| + |C| \geq 3$ and,
without loss of generality, $|B| \geq 2$. Let $b \in B$. Applying
the submodularity of $r$ twice, we get that $d = r(A \cup B \cup
C) - r(A \cup B) \geq r(A \cup b \cup C) - r(A \cup b) \geq r(A
\cup C) - r(A) = d.$ It follows that $r(A \cup b \cup C) - r(A
\cup b) = d$ also.

We can apply the induction hypothesis to the sets $A, \{b\},C$,
since $A \cup b, A \cup C \in \C$ and $|\{b\}| + |C| < |B| + |C|$.
We conclude that $A \cup b \cup C \in \C$. We can then apply the
induction hypothesis to the sets $A \cup b, B-b, C$, since $A
\cup B, A \cup b \cup C \in \C$ and $|B-b| + |C| < |B| + |C|$. We
conclude that $A \cup B \cup C \in \C$, as desired. $\Box$

\medskip

Propositions \ref{prop:semi-mat} and \ref{prop:modideal} show us
how to obtain a pair $(M,\I)$ of a matroid $M$ and one of its
modular ideals $\I$, given a semimatroid $\C$. Now we show that it
is possible to recover $\C$ from the pair $(M, \I)$.

\begin{proposition} \label{prop:mat,modid-semimat}
Let $M=(S,r)$ be a matroid, and let $\I$ be a modular ideal of
$M$. Let $r_{\I}$ be the restriction of the rank function of $M$
to $\I$. Then $(S, \I, r_{\I})$ is a semimatroid.
\end{proposition}

\noindent \emph{Proof.} The rank function $r_{\I}$ inherits axioms
(R1)-(R3) from $r_M$. (CR1) is easy. If $X,Y \in \I$ and
$r_{\I}(X) = r_{\I}(X \cap Y)$, then $r(Y) = r(X \cup Y)$ by
submodularity. Thus $\{X,Y\}$ is a modular pair in $M$, and $X
\cup Y \in \I$.

Now we check (CR2). We start by showing that $\I$ must contain
every independent set of $M$. In fact, assume that $I$ is a
minimal independent set which is not in $\I$. Since $\I$ contains
all non-loop elements, $I$ has at least two elements $a$ and $b$.
Then $\I$ contains the modular pair $\{I-a, I-b\}$, so it
contains their union $I$, a contradiction.

Now take $X,Y \in \I$ with $|X| < |Y|$, and pick $y \in Y$ such
that $r(X \cup y) = r(X) + 1$. Let $X'$ be an independent subset
of $X$ of rank $r(X)$; then $X' \cup y$ is an independent set of
rank $r(X)+1$. Therefore $\I$ contains the modular pair $\{X' \cup
y,X\}$, so it contains their union $X \cup y$. $\Box$

\medskip

\begin{theorem} \label{th:semimatmatid}
Let $S$ be a finite set. Let $\semimat(S)$ be the set of
semimatroids on $S$. Let $\matid(S)$ be the set of pairs $(M,\I)$
of a matroid $M$ on $S$ and a modular ideal $\I$ of $M$.
\begin{enumerate}
\item
The assignment $(S, \C, r_{\C}) \mapsto (M_{\C}, \C)$ is a map
$\semimat(S) \rightarrow \matid(S)$.
\item
The assignment $(M, \I) \mapsto (S, \I, r_{\I})$ is a map
$\matid(S) \rightarrow \semimat(S)$.
\item
The two maps above are inverses, and give a one-to-one
correspondence between $\semimat(S)$ and $\matid(S)$.
\end{enumerate}
\end{theorem}

\noindent \emph{Proof.} The first and second parts are
restatements of Propositions \ref{prop:semi-mat} and
\ref{prop:modideal} and Proposition \ref{prop:mat,modid-semimat},
respectively.

Denote the maps $\semimat(S) \rightarrow \matid(S)$ and $\matid(S)
\rightarrow \semimat(S)$ above by $f$ and $g$ respectively. It is
immediate that $g \circ f$ is the identity map in $\semimat(S)$.
To check that $f \circ g$ is the identity map in $\matid(S)$, we
need to show the following. Given a matroid $M = (S,r)$ and a
modular ideal $\I$ of $M$, $r(X) = \max\{r(Y) \, | \, Y \subseteq
X \, , \, Y \in \I\}$ for all $X \subseteq S$. But this is easy:
it is clear that $r(X) \geq \max\{r(Y) \, | \, Y \subseteq X \, ,
\, Y \in \I\}$. Equality is attained because $X$ has an
independent subset $X'$ of rank $r(X)$; since $X'$ is
independent, it is in $\I$. $\Box$

\medskip

Before we continue our analysis, we state explicitly a simple
property of semimatroids and modular ideals which is implicit in
the proofs above.

In a semimatroid $(S, \C, r_{\C})$, all the maximal sets in $\C$
have the same rank, which we denote $r_{\C}$. In a modular ideal
$\I$ of a matroid $M=(S,r)$, all the maximal sets have maximum
rank $r=r(S)$.

\section{Elementary preimages and single-element coextensions.}
\label{sec:preim}

Now we show that a semimatroid is also equivalent to a pair
$(M,M')$ of a matroid $M$ and one of its rank-increasing
single-element coextensions $M'$.
To do it, we outline the correspondence between the modular
ideals, the elementary preimages and the rank-increasing
single-element coextensions of a matroid.

This correspondence is just the dual of the well understood
correspondence between the modular filters, the elementary
quotients, and the rank-preserving single-element extensions of a
matroid \cite{Cr65}, \cite{Do76}, \cite{Ku86}. Therefore we omit
all the proofs of these results, and refer the reader to the
relevant literature.

\begin{defin}
A \emph{quotient map} $N \rightarrow M$ is a pair of matroids $M,
N$ on the same ground set such that every flat of $M$ is a flat
of $N$.
\end{defin}

There are several other equivalent definitions of quotient maps,
including the following.

\begin{proposition} \cite[Proposition 8.1.6]{Ku86}
Let $M$ and $N$ be two matroids on the set $S$. The following are
equivalent: \\
(i) $N \rightarrow M$ is a quotient map. \\
(ii) For any $A \subseteq S$, $\cl_N(A) \subseteq
\cl_M(A)$. \\
(iii) For any $A \subseteq B \subseteq S$, $r_N(B) - r_N(A) \geq
r_M(B) - r_M(A)$.
\end{proposition}

\begin{defin}
An \emph{elementary quotient map} is a quotient map $N \rightarrow
M$ such that $r(N)-r(M) = 0$ or $1$.
\end{defin}

We will focus our attention on elementary quotient maps. Their
importance is the following. Perhaps the most useful notion of a
morphism in the category of matroids is that of a \emph{strong
map}. Every strong map between matroids can be regarded
essentially as a quotient map, followed by an embedding of a
submatroid into a matroid. Also, every quotient map can be
factored into a sequence of elementary quotient maps. Therefore,
elementary quotient maps are essentially the building blocks of
strong maps. For more information on this topic, we refer the
reader to \cite{Ku86}.

\begin{defin}
An \emph{elementary preimage} of a matroid $M$ is a matroid $N$
on the same ground set such that $N \rightarrow M$ is an
elementary quotient map.
\end{defin}

The following proposition explains the relevance of elementary
preimages and quotient maps in our investigation.

\begin{theorem} \cite[Proposition 6.5]{Do76} \label{th:idealpreim}
Let $M = (S,r_M)$ be a  matroid. Let $\ideal(M)$ be the set of
modular ideals of $M$ and let $\preim(M)$ be the set of elementary
preimages of $M$.
\begin{enumerate}
\item
Given $\I \in \ideal(M)$, define the rank function $r_N:2^S
\rightarrow \NN$ by:
\begin{equation*}
r_N(A) = \left\{ \begin{array}{ll}
                r_M(A) & \mbox{if } A \in \I \\
                r_M(A)+1 & \mbox{if } A \notin \I
                \end{array} \right.
\end{equation*}
Then $N = (S, r_N)$ is a matroid, and $N \in \preim(M)$.

\item
Given $N \in \preim(M)$, let $\I= \{A \in S: r_M(A) = r_N(A)\}.$
Then $\I \in \ideal(M)$.

\item
The two maps $\ideal(M) \rightarrow \preim(M)$ and $\preim(M)
\rightarrow \ideal(M)$ defined above are inverses. They establish
a one-to-one correspondence between $\ideal(M)$ and $\preim(M)$.
\end{enumerate}
\end{theorem}

\begin{corollary} \label{th:semimatmatpreim}
Given a finite set $S$, let $\matpreim(S)$ be the set of pairs
$(M,N)$ of a matroid $M$ on $S$ and one of its elementary
preimages $N$. Then there are one-to-one correspondences between
$\semimat(S)$, $\matid(S)$ and $\matpreim(S)$.
\end{corollary}

\noindent \emph{Proof.} Combine Theorems \ref{th:semimatmatid} and
\ref{th:idealpreim}. $\Box$.

\medskip

\begin{defin}
Let $M$ be a matroid on the ground set $S$ and let $p$ be an
element not in $S$. A \emph{single-element coextension} of $M$ by
$p$ is a matroid $\Nh$ on the set $S \cup p$ such that $M =
\Nh/p$. $\Nh$ is \emph{rank-increasing} if $r(\Nh) > r(M)$.
\end{defin}

It is worth remarking that most single-element coextensions of
$M$ by $p$ are rank-increasing. The only one which is not
rank-increasing is the matroid $\Nh$ on $S \cup p$ such that
$r_{\Nh}(A \cup p) = r_{\Nh}(A) = r_M(A)$ for all $A \subseteq
S$; \emph{i.e.}, the one where $p$ is a loop.

\begin{theorem} \cite[dual of Theorem 8.3.2]{Ku86} \label{th:preimcoext}
Let $M$ be a matroid and $p$ be an element not in its ground set.
Let $\coext(M,p)$ be the set of rank-increasing single-element
coextensions of $M$ by $p$.
\begin{enumerate}
\item
Given $N \in \preim(M)$, define $r_{\Nh}:2^{S \cup p} \rightarrow
\NN$ by
\begin{eqnarray*}
r_{\Nh}(A) &=& r_N(A)\\
r_{\Nh}(A \cup p) &=& r_M(A)+1
\end{eqnarray*}
for $A \subseteq S$. Then $\Nh = (S \cup p, r_{\Nh})$ is a
matroid, and $\Nh \in \coext(M,p)$.
\item
If $\Nh \in \coext(M,p)$, then the matroid $N = \Nh - p$ is in
$\preim(M)$.
\item
The two maps $\preim(M) \rightarrow \coext(M,p)$ and $\coext(M,p)
\rightarrow \preim(M)$ defined above are inverses. They establish
a one-to-one correspondence between $\preim(M)$ and $\coext(M,p)$.
\end{enumerate}
\end{theorem}

\begin{corollary} \label{th:semimatmatcoext}
Given a finite set $S$ and an element $p \notin S$, let
$\matcoext(S,p)$ be the set of pairs $(M,\Nh)$, where $M$ is a
matroid on $S$ and $\Nh$ is one of its rank-increasing
single-element coextensions by $p$. Then there are one-to-one
correspondences between $\semimat(S)$, $\matid(S)$,
$\matpreim(S)$ and $\matcoext(S,p)$.
\end{corollary}

\noindent \emph{Proof.} Combine Theorems \ref{th:semimatmatid},
\ref{th:idealpreim} and \ref{th:preimcoext}. $\Box$

\medskip

We briefly mention that given a matroid $M$, there are other
objects in correspondence with the modular ideals of $M$. Two
such examples are the \emph{modular cocuts} of $M$ and the
\emph{colinear subclasses} of $M$. They are the duals of modular
cuts and linear subclasses, respectively.

A modular cocut $\U$ of $M$ is a collection of circuit unions of
$M$ satisfying two conditions. First, if $U_1 \subseteq U_2$ are
circuit unions and $U_2 \in \U$, then $U_1 \in \U$. Second, if
$U_1, U_2 \in \U$ and $\{U_1, U_2\}$ is a modular pair in $M$,
then $U_1 \cup U_2 \in \U$.

A colinear subclass $\C$ of $M$ is a set of circuits of $M$ such
that if $C_1,C_2 \in \C$ and $r(C_1 \cup C_2) = |C_1 \cup C_2| -
2$, and $C_3 \subseteq C_1 \cup C_2$ is a circuit, then $C_3 \in
\C$.

The details and proofs of the (dual) correspondences appear in
\cite[Theorem 7.2.2]{Ox92} and \cite{Cr65}, respectively.

\medskip

\begin{figure}[!h]
\centering
\includegraphics[height=5in]{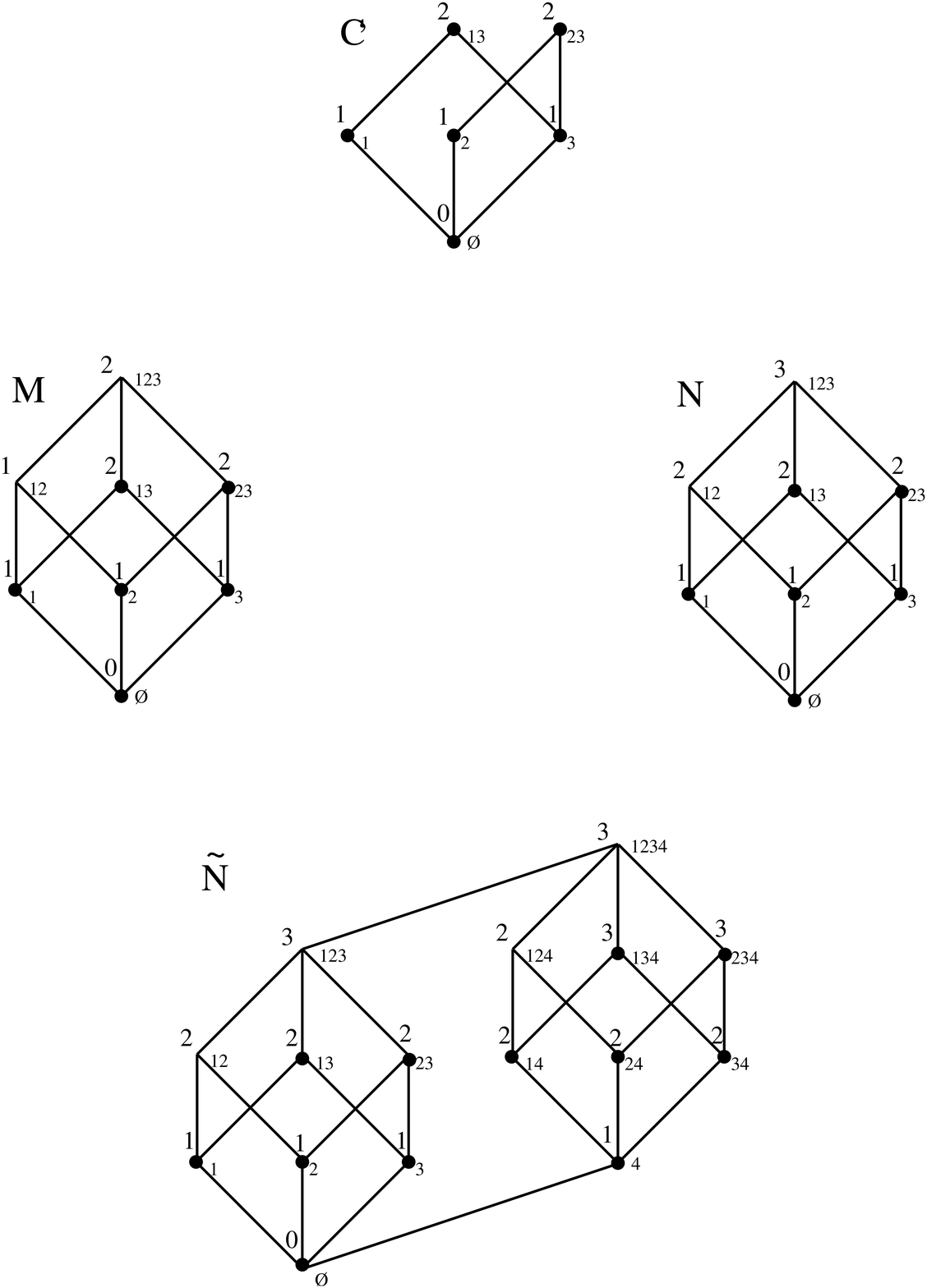}
\caption{The semimatroid $\C$ and its corresponding matroids.}
\label{fig:bijections}
\end{figure}

We end this section by reviewing the correspondences and objects
of Sections \ref{sec:modid} and \ref{sec:preim} with an example.
Let $\C=(S, \C, r)$ be the semimatroid consisting of the set
$S=[3]$, the collection of central sets $\C = 2^{[3]} - \{12,
123\}$, and the rank function $r(X) = |X|$ for $X \in \C$. It is
easy to check that this is, indeed, a semimatroid.  The first
diagram of Figure \ref{fig:bijections} depicts the poset on $\C$
ordered by inclusion; below and to the right of each node we have
written the set in $\C$ corresponding to it, and above and to its
left we have written its rank.

To the semimatroid $\C$, we have assigned a pair $(M, \I) \in
\matid(S)$, a pair $(M,N) \in \matpreim(S)$ and a pair $(M, \Nh)
\in \matcoext(S)$. To obtain the matroid $M$, we add the subsets
of $[3]$ not in $\C$ to the diagram above, to get the Boolean
algebra $2^{[3]}$. We have placed big nodes on the sets of this
poset which are in the diagram of $\C$, and small nodes on the new
sets. To obtain the rank function of $M$, we copy the rank
function of $\C$ on the big nodes. On each small node, we put the
largest number that we can find on a big node below it. The big
nodes form the modular ideal $\I$ of $M$.

To obtain the matroid $N$, we simply leave the rank function of
$M$ fixed on the big nodes, and increase it by $1$ on the little
nodes.

Finally, to obtain the matroid $\Nh$, we glue two Boolean algebras
$2^{[3]}$, to obtain a  Boolean algebra $2^{[4]}$ on $4$ elements.
(We have omitted most of the  ``diagonal'' edges of this poset for
clarity.) On the lower copy of the Boolean algebra, we put the rank
function of $N$. On the upper copy, we put the rank function of $M$,
increased by $1$.

\section{Pointed matroids.} \label{sec:pointed}

We now establish a correspondence between semimatroids and
pointed matroids.

\begin{defin} \cite{Br71}
A \emph{pointed matroid} is a pair $(M,p)$ of a matroid $M$ and a
distinguished element $p$ of its ground set.
\end{defin}

Pointed matroids are a combinatorial tool often used in the study
of affine hyperplane arrangements. The connection between them is
the following. Consider an affine arrangement $\A = \{H_1, \ldots,
H_k\}$ in $\RR^n$, where $H_i$ is defined by the equation $v_i
\cdot x = c_i$.

\begin{defin} \label{def:cone}
The \emph{cone over $\A$} is the arrangement $c\A = \{H_1',
\ldots, H_k', H\}$ in $\RR^{n+1}$, where $H_i'$ is defined
\footnote{with a slight abuse of notation} by the equation $v_i
\cdot x = c_ix_{n+1}$ for $1 \leq i \leq k$, and $H$ is the \emph{
additional hyperplane} $x_{n+1}=0$.
\end{defin}

Being a central arrangement, $c\A$ has a matroid $M_{c\A}$ on the
ground set $c\A$ associated to it. To the arrangement $\A$, we
associate the pointed matroid $(M_{c\A},H)$.

\begin{theorem} \label{th:matcoextpointedmat}
Let $S$ be a set and let $p \notin S$. Let $\pointedmat(S,p)$ be
the set of pointed matroids $(M,p)$ on $S \cup p$ such that $p$ is
not a loop of $M$. There are one-to-one correspondences between
$\semimat(S)$, $\matid(S)$, $\matpreim(S)$, $\matcoext(S,p)$ and
$\pointedmat(S,p)$.
\end{theorem}

\noindent \emph{Proof.} It suffices to show a correspondence
between $\matcoext(S,p)$ and $\pointedmat(S,p)$. The elements of
$\matcoext(S,p)$ are the pairs $(\Nh/p, \Nh)$ for all matroids
$\Nh$ on $S \cup p$ such that $r(\Nh) > r(\Nh/p)$; \emph{i.e.},
such that $p$ is not a loop. The map $(\Nh/p, \Nh) \mapsto
(\Nh,p)$ establishes the desired bijection. $\Box$

\medskip

At this point, given a set $S$ and an element $p \notin S$, we
have bijections between $\semimat(S)$, $\matid(S)$,
$\matpreim(S)$, $\matcoext(S,p)$ and $\pointedmat(S,p)$, provided
by Theorems \ref{th:semimatmatid}, \ref{th:semimatmatpreim},
\ref{th:semimatmatcoext}, and \ref{th:matcoextpointedmat}. The
bijection $\pointedmat(S,p) \rightarrow \semimat(S)$ is an
important one. We have obtained it as the composition of four
bijections, and now we wish to describe it explicitly.

\begin{theorem}\label{th:pointedmatsemimat}
Let $S$ be a set and let $p \notin S$.
\begin{enumerate}
\item
For $(\Nh,p) \in \pointedmat(S,p)$, let $\C = \{A \subseteq S \, |
\, p \notin \cl_{\Nh}(A)\}$ and let $r_{\C}$ be the restriction of
$r_{\Nh}$ to $\C$. Then $(S, \C, r_{\C})$ is a semimatroid.
\item
For $(S, \C, r_{\C}) \in \semimat(S)$, define $r_{\Nh}:2^{S \cup
p} \rightarrow \NN$ by
\begin{eqnarray*}
r_{\Nh}(A) &=& \left\{ \begin{array}{ll}
                r_{\C}(A) & \mbox{if } A \in \C \\
                \max \{r_{\C}(B) \, | \, B \subseteq A, B \in \C\} + 1 & \mbox{if } A \notin \C
                \end{array} \right. \\
r_{\Nh}(A \cup p) &=& \left\{ \begin{array}{ll}
                r_{\Nh}(A)+1 & \mbox{if } A \in \C\\
                r_{\Nh}(A) & \mbox{if } A \notin \C
                \end{array} \right.
\end{eqnarray*}
for $A \subseteq S$. Then $r_{\Nh}$ is a rank function on $S \cup
p$, and $(\Nh,p) \in \pointedmat(S,p)$.
\item
The two maps $\pointedmat(S,p) \rightarrow \semimat(S)$ and
$\semimat(S) \rightarrow \pointedmat(S,p)$ defined above are
inverses. They establish a one-to-one correspondence between
$\pointedmat(S,p)$ and $\semimat(S)$.
\end{enumerate}
\end{theorem}

\noindent \emph{Proof.} We will show that, if we start with
$(\Nh,p) \in \pointedmat(S,p)$ and trace the bijections of
Theorems \ref{th:matcoextpointedmat}, \ref{th:preimcoext},
\ref{th:idealpreim} and \ref{th:semimatmatid}, we obtain the
semimatroid $\C(\Nh,p)$.

Under the bijection of Theorem \ref{th:matcoextpointedmat},
$(\Nh,p) \in \pointedmat(S,p)$ corresponds to $(\Nh/p,\Nh) \in
\matcoext(S,p)$.

Under the bijection of Theorem \ref{th:preimcoext}, $\Nh \in
\coext(\Nh/p)$ corresponds to $\Nh-p \in \preim(\Nh/p)$.

$\Nh-p \in \preim(\Nh/p)$, under the bijection of Theorem
\ref{th:idealpreim}, corresponds to the modular ideal $\C = \{A
\subseteq S \, | \, r_{\Nh/p}(A) = r_{\Nh-p}(A)\} \in
\ideal(\Nh/p)$. Since $p$ is not a loop of $\Nh$, $r_{\Nh/p}(A) =
r_{\Nh}(A \cup p) -1$ and $r_{\Nh-p}(A) = r_{\Nh}(A)$. Therefore
$\C = \{A \subseteq S \, | \, p \notin \cl_{\Nh}(A)\}$.

Finally, under the bijection of Theorem \ref{th:semimatmatid},
$(\Nh/p, \C) \in \matid(S)$ corresponds to $(S, \C, r_{\C}) \in
\semimat(S)$.

Similarly, if we start with a semimatroid $(S, \C, r_{\C})$ and
keep track of its successive images under the bijections of
Theorems \ref{th:semimatmatid}, \ref{th:idealpreim},
\ref{th:preimcoext} and \ref{th:matcoextpointedmat}, we get the
pointed matroid $(\Nh,p)$ described.

This theorem then becomes a consequence of the previous ones.
$\Box$

\medskip

It is not difficult to see that, under the coning construction,
the central subsets of a hyperplane arrangement $\A$ correspond
to the subsets of $c\A$ whose closure in $M_{c\A}$ does not
contain the additional hyperplane $H$. Theorem
\ref{th:pointedmatsemimat} shows that, for semimatroids, the
natural analog of the cone of a semimatroid $\C$ is the matroid
$\Nh$ of the pointed matroid $(\Nh,p)$ corresponding to it.

\medskip

The triple of matroids $(\Nh, \Nh-p, \Nh/p) = (\Nh, N, M)$ is
sometimes called the \emph{triple of the pointed matroid}
$(\Nh,p)$. We will also call it the \emph{triple of the
semimatroid} $\C$.

%
%
%

\section{Geometric semilattices.} \label{sec:semilat}

We now discuss geometric semilattices and their relationship to
semimatroids. We start by recalling some poset terminology. For
more background information, see for example \cite[Chapter
3]{St86}.

A \emph{meet semilattice} is a poset $K$ such that any subset $S
\subseteq K$ has a \emph{greatest lower bound} or \emph{meet}
$\wedge S$: an element such that $\wedge S \leq s$ for all $s \in
S$, and $\wedge S \geq t$ for any $t \in K$ such that $t \leq s$
for all $s \in S$. Such posets have a minimum element $\hat{0}$.

Notice that if a set $S$ of elements of a meet semilattice has an
upper bound, then it has a least upper bound, or \emph{join}
$\vee S$. It is the meet of the upper bounds of $S$.

A \emph{lattice} is a poset $L$ such that any subset $S \subseteq
L$ has a greatest lower bound and a least upper bound. Clearly, if
a meet semilattice has a maximum element, then it is a lattice.

A meet semilattice $K$ is \emph{ranked} with \emph{rank function}
$r:K \rightarrow \NN$ if, for all $x \in K$, every maximal chain
from $\hat{0}$ to $x$ has the same length $r(x)$. An \emph{atom}
is an element of rank $1$. A set of atoms $A$ is
\emph{independent} if it has an upper bound and $r(\vee A) = |A|$.

\begin{defin}
A \emph{geometric semilattice} is a ranked meet semilattice
satisfying the following two conditions.

\noindent {\bf(G1)} Every element is a join of atoms.

\noindent {\bf(G2)} The collection of independent set of atoms is
the collection of independent sets of a matroid.

A \emph{geometric lattice} is a ranked lattice satisfying (G1) and
(G2).
\end{defin}

Geometric lattices arise very naturally in matroid theory from
the following result. Recall that a matroid $M = (S,r)$ is
\emph{simple} if $r(x)=1$ for all $x \in S$ and $r(\{x,y\})=2$
for all $x,y \in S, x \neq y$.

\begin{theorem} \cite{Bi35}, \cite{Cr70} \label{th:latmat}
A poset is a geometric lattice if and only if it is isomorphic to
the poset of flats of a matroid. Furthermore, each geometric
lattice is the poset of flats of a unique simple matroid, up to
isomorphism.
\end{theorem}

From this point of view, semimatroids are the ``right"
generalization of matroids, as the following theorem shows.

\begin{defin}
A semimatroid $\C = (S,\C,r_{\C})$ is \emph{simple} if $\{x\} \in
\C$ and $r_{\C}(x)=1$ for all $x \in S$, and $r_{\C}(\{x,y\})=2$
for all $\{x,y\} \in \C$ with $x \neq y$.
\end{defin}

\begin{theorem} \label{th:semilatsemimat}
A poset is a geometric semilattice if and only if it is
isomorphic to the poset of flats of a semimatroid. Furthermore,
each geometric semilattice is the poset of flats of a unique
simple semimatroid, up to isomorphism.
\end{theorem}

To prove Theorem \ref{th:semilatsemimat} we use the following two
propositions.

\begin{proposition} \label{prop:semilatlat}\cite{Wa86}
A poset $K$ is a geometric semilattice if and only if there is a
geometric lattice $L$ with an atom $p$ such that $K = L - [p,
\hat{1}]$. \footnote{Here $[p, \hat{1}]$ denotes the interval of
elements greater than or equal to $p$ in the poset $L$.}
Furthermore, $L$ and $p$ are uniquely determined by $K$.
\end{proposition}

\begin{proposition} \label{prop:flatsofsemimat}
Let $\C=(S, \C, r_{\C})$ be a semimatroid, and let $(\Nh,p)$ be
the pointed matroid on $S \cup p$ corresponding to it under the
bijection of Theorem \ref{th:pointedmatsemimat}. Let $K(\C)$ and
$L(\Nh)$ be the posets of flats of $\C$ and $\Nh$. Then $K(\C) =
L(\Nh) - [p, \hat{1}]$.
\end{proposition}

\noindent \emph{Proof.} Since both posets are ordered by
containment, we only need to show the equality of the sets
$K(\C)$ and $L(\Nh) - [p, \hat{1}]$.

First we show that $K(\C) \subseteq L(\Nh) - [p, \hat{1}]$. Let $X
\in K(\C)$. Then for all $x \notin X$ such that $X \cup x \in
\C$, $r_{\C}(X \cup x) = r_{\C}(X) + 1$, and therefore $r_{\Nh}(X
\cup x) = r_{\Nh}(X) + 1$. To check that $X$ is a flat in $\Nh$,
we need to show that this equality still holds if $X \cup x
\notin \C$. This is not difficult: if that is the case and $x
\neq p$, then $r_{\Nh}(X \cup x) = \max \{r_{\C}(Y) \, | \, Y
\subseteq X \cup x, Y \in \C\} + 1 \geq r_{\C}(X)+1 =
r_{\Nh}(X)+1$. Clearly then equality must hold. The case $x = p$
is easier, but needs to be checked separately.

Hence $K(\C) \subseteq L(\Nh)$, and since no element of $\C$
contains $p$, $K(\C) \subseteq L(\Nh) - [p, \hat{1}]$.

The inverse inclusion is easier. If $X$ is a flat in $\Nh$ not
containing $p$, then $r_{\Nh}(X \cup x) = r_{\Nh}(X)+1$ for all
$x \notin X$. When $X \cup x \in \C$, this equality says that
$r_{\C}(X \cup x) = r_{\C}(X)+1$. Therefore $X$ is a flat in $\C$
also. $\Box$

\medskip

\noindent \emph{Proof of Theorem \ref{th:semilatsemimat}.} It is not
difficult to check that the bijection of Theorem
\ref{th:pointedmatsemimat} restricts to a bijection between
\emph{simple pointed matroids} (pointed matroids $(\Nh,p) \in
\pointedmat(S,p)$ such that $\Nh$ is simple) and simple semimatroids.
The result then follows combining this fact with
Theorem \ref{th:latmat} and Propositions
\ref{prop:semilatlat} and \ref{prop:flatsofsemimat}. $\Box$

\section{Duality, deletion and contraction.} \label{sec:del-cont}

Like matroids, semimatroids have natural notions of duality,
deletion and contraction, which we now define.

\begin{defin}
Let $\C = (S, \C, r_{\C})$ be a semimatroid. Extend the function
$r_{\C}$ to a matroid rank function $r:2^S \rightarrow \NN$ as in
Proposition \ref{prop:semi-mat}. Define the simplicial complex
$\C^* = \{X \subseteq S \, | \, S-X \notin \C\}$, and the rank
function $r^*:\C^* \rightarrow \NN$ by $r^*(X) = |X| - r +
r(S-X)$. The \emph{dual of} $\C$ is the triple $\C^* = (S, \C^*,
r^*)$.
\end{defin}

\begin{proposition}
$\C^*$ is a semimatroid.
\end{proposition}

\noindent \emph{Proof.} It is possible to simply check that
$\C^*$ satisfies the axioms of a semimatroid. It is shorter to
proceed as follows.

Consider the pair $(M,N) \in \matpreim(S)$ associated to $\C$
under Corollary \ref{th:semimatmatpreim}. It is known
\cite[Proposition 8.1.6(f)]{Ku86} that if $N$ is an elementary
preimage of $M$, then $M^*$ is an elementary preimage of $N^*$.
From the pair $(N^*,M^*) \in \matpreim(S)$, we then get a
semimatroid using Corollary \ref{th:semimatmatpreim} again. It is
straightforward to check that this semimatroid is precisely
$\C^*$. $\Box$

\medskip

\begin{proposition}
For any semimatroid $\C$, we have that $(\C^*)^* = \C$.
\end{proposition}

\noindent \emph{Proof.} This is easy to check directly from the
definition. $\Box$

\medskip

\begin{defin}
Let $\C = (S, \C, r_{\C})$ be a semimatroid and let $e \in S$ be
such that $\{e\} \in \C$. Let $\C/e = \{A \subseteq S-e \, | \, A
\cup e \in \C\}$ and, for $A \in \C/e$, let $r_{\C/e}(A) =
r_{\C}(A \cup e) - r_{\C}(e)$. The \emph{contraction of $e$ from
$\C$} is the triple $\C/e = (S-e, \C/e, r_{\C/e})$.
\end{defin}
\begin{proposition}
$\C/e$ is a semimatroid.
\end{proposition}

\noindent \emph{Proof.} Checking the axioms of a semimatroid is
straightforward. $\Box$

\medskip

\begin{defin}
Let $\C = (S, \C, r_{\C})$ be a semimatroid and let $e \in S$ be
such that $\{e\} \in \C$. Let $\C-e = \{A \in \C \, | \, e \notin
A\}$ and, for $A \in \C-e$, let $r_{\C-e}(A) = r_{\C}(A)$. The
\emph{deletion of $e$ from $\C$} is the triple $\C-e = (S-e,
\C-e, r_{\C-e})$.
\end{defin}
\begin{proposition}
$\C-e$ is a semimatroid.
\end{proposition}

\noindent \emph{Proof.} Checking the axioms of a semimatroid is
straightforward. $\Box$

\medskip

Again, as with matroids, there are two special kinds of elements
that we need to pay special attention to when we perform deletion
and contraction.

\begin{defin}
A \emph{loop} of a semimatroid $\C = (S, \C, r_{\C})$ is an
element $e \in S$ such that $\{e\} \in \C$ and $r_{\C}(e) = 0$.
\end{defin}

\begin{defin}
An \emph{isthmus} of a semimatroid $\C = (S, \C, r_{\C})$ is an
element $e \in S$ such that, for all $A \in \C$, $A \cup e \in
\C$ and $r_{\C}(A \cup e) = r_{\C}(A) + 1$.
\end{defin}

\begin{lemma} \label{lemma:loop}
If $e \in S$ is a loop of the semimatroid $\C = (S, \C, r_{\C})$,
then $r_{\C/e} = r_{\C}$. Otherwise, $r_{\C/e} = r_{\C} - 1$.
\end{lemma}

\noindent \emph{Proof.} Clearly $r_{\C/e} \leq r_{\C}$. If $e$ is
a loop, consider any $A \in \C$. (CR1') applies to $\{e\}$ and
$A$, so $A \cup e \in \C$ and $r_{\C}(A \cup e) = r_{\C}(A)$.
Therefore the maximum rank $r_{\C}$ in $\C$ is achieved for some
$A \in \C/e$. But then we have $r_{\C/e}(A) = r_{\C}(A \cup e) - 0
= r_{\C}$,
so $r_{\C/e} = r_{\C}$.

If $e$ is not a loop, then for all $A \in \C/e$ we have
$r_{\C/e}(A) = r_{\C}(A \cup e) - 1$, so $r_{C/e} \leq r_C - 1$.
Equality holds: if we start with $\{e\} \in \C$ and repeatedly
apply (CR2') with an element of $\C$ of rank $r_{\C}$, we can
obtain a set $A \cup e$ of rank $r_{\C}$. Then $r_{C/e}(A) = r_C -
1$. $\Box$

\medskip

\begin{lemma} \label{lemma:isthmus}
If $e \in S$ is an isthmus of the semimatroid $\C = (S, \C,
r_{\C})$, then $r_{\C-e} = r_{\C}-1$. Otherwise, $r_{\C-e} =
r_{\C}$.
\end{lemma}

\noindent \emph{Proof.} Clearly $r_{\C-e} \leq r_{\C}$. If $e$ is
an isthmus then it is clear from the definition that $r_{\C-e} =
r_{\C}-1$.

If $e$ is not an isthmus, there are two cases. If there is an $A
\in \C$ such that $A \cup e \notin \C$, take a maximal one. It is
also a maximal set in $\C$, so it has maximum rank $r_{\C}$; and $A \in
\C-e$, so $r_{\C-e} = r_{\C}$. The other possibility is that for
all $A \in \C$, we have $A \cup e \in \C$ and $r(A \cup e) =
r(A)$. In this case it is also clear that $r_{\C-e} = r_{\C}$.
$\Box$

\medskip

\begin{lemma} \label{lemma:loopisthmus}
If $e \in S$ is a loop or an isthmus of the semimatroid $\C = (S,
\C, r_{\C})$, then $\C-e = \C/e$.
\end{lemma}

\noindent \emph{Proof.} This is clear from Lemmas
\ref{lemma:loop} and \ref{lemma:isthmus} and their proofs. $\Box$

\medskip

\section{The Tutte polynomial.} \label{sec:Tuttesemi}

With the background results that we have established, we are now
able to define and study the Tutte polynomial of a semimatroid.

\begin{defin}
The \emph{Tutte polynomial} of a semimatroid $\C= (S, \C, r_{\C})$
is defined by
\begin{equation}\label{eq:Tutte}
T_{\C}(x,y) = \sum_{X \, \in \, \C}
(x-1)^{r_{\C}-r_{\C}(X)}(y-1)^{|X|-r_{\C}(X)}.
\end{equation}
\end{defin}

If $\A$ be a hyperplane arrangement and $\C_{\A}$ is the
semimatroid determined by it, then the Tutte polynomial of the
semimatroid $\C_{\A}$ is precisely the Tutte polynomial of the
arrangement $\A$, as defined and studied in \cite{arrangements}.
That paper focuses on enumerative aspects arising from the
computation of these polynomials; here we will concentrate our
attention on matroid-theoretical considerations.

\medskip

\begin{figure}[b]
\centering
\includegraphics[width=3.5in]{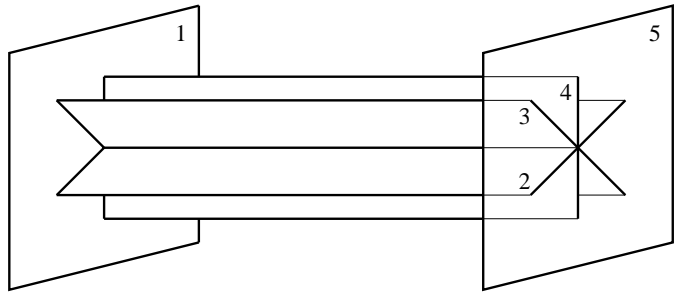}
\caption{The arrangement $\A$.} \label{fig:A}
\end{figure}

\noindent \emph{Example.} Figure \ref{fig:A} shows a hyperplane
arrangement $\A$ in $\RR^3$, consisting of the five planes
$x+y+z=0, x=y, y=z, z=x$ and $x+y+z=1$ in that order.

Table \ref{table:Tutte} shows all the central subsets of $\A$, and
their contributions to the Tutte polynomial of $\A$.
\begin{table}[t]
\centering \setlength{\extrarowheight}{2pt}
\begin{tabular}{|r|l|}
\hline
central subset of $\A$ &contribution to $T_{\A}(x,y)$ \\
\hline $\emptyset$  & $(x-1)^3(y-1)^0$ \\
\hline
$1,2,3,4,5$ & $(x-1)^2(y-1)^0$ \\
\hline
$12,13,14,23,24,25,34,35,45$  & $(x-1)^1(y-1)^0$ \\
\hline
$123, 124, 134, 235, 245, 345$ & $(x-1)^0(y-1)^0$ \\
\hline
$234$ & $(x-1)^1 (y-1)^1$ \\
\hline
$1234, 2345$ & $(x-1)^0 (y-1)^1$ \\
\hline
\end{tabular}
\caption{Computing the Tutte polynomial $T_{\A}(x,y).$}
\label{table:Tutte}
\end{table}
We find that
\begin{eqnarray*}
T_{\A}(x,y) &=& (x-1)^3 + 5(x-1)^2 + 9(x-1) + 6 + (x-1)(y-1) + 2(y-1) \\
& = & x^3 + 2x^2 + x y + x + y.
\end{eqnarray*}

As in the matroid setting, the Tutte polynomial of a semimatroid
satisfies the following simple recursive formula, known as the
\emph{deletion-contraction} relation.

\begin{proposition}\label{pr:del-cont}
Let $\C = (S, \C, r_{\C})$ be a semimatroid, and let $e \in S$ be
such that $\{e\} \in \C$.
\begin{enumerate}
\item[(i)]
$T_{\C}(x,y) = T_{\C-e}(x,y) + T_{\C/e}(x,y)$ if $e$ is neither an
isthmus nor a loop and $\{e\} \in \C$.
\item[(ii)]
$T_{\C}(x,y) = x\,T_{\C-e}(x,y)$ if $e$ is an isthmus.
\item[(iii)]
$T_{\C}(x,y) = y\,T_{\C/e}(x,y)$ if $e$ is a loop.
\item[(iv)]
If $e \in S$ and $\{e\} \notin \C$ then $T_{(S, \C, r_{\C})}(x,y)
= T_{(S-e, \C, r_{\C})}(x,y)$.
\end{enumerate}
\end{proposition}

\noindent \emph{Proof.} We have

\begin{eqnarray*}
T_{\C}(x,y) &=& \sum_{\stackrel{X \, \in \, \C}{e \notin X}}
(x-1)^{r_{\C}-r_{\C}(X)}(y-1)^{|X|-r_{\C}(X)} + \\
&& \hspace{2cm} \sum_{X \cup e \, \in \, \C}
(x-1)^{r_{\C}-r_{\C}(X \cup e)}(y-1)^{|X \cup e|-r_{\C}(X \cup
e)}.
\end{eqnarray*}

Notice that, if $r_{\C} = r_{\C-e}$, the first sum in the right
hand side is exactly the Tutte polynomial of $\C-e$. If, on the
other hand, $r_{\C} = r_{\C-e}+1$, the only difference is that we
get an extra factor of $(x-1)$. More precisely, in view of Lemma
\ref{lemma:isthmus}, the first sum of the right hand side is
$T_{\C-e}(x,y)$ if $e$ is not an isthmus, and
$(x-1)T_{\C-e}(x,y)$ if it is an isthmus.

Similarly, from Lemma \ref{lemma:loop}, the second sum is
$T_{\C/e}(x,y)$ if $e$ is not a loop, and $(y-1)T_{\C/e}(x,y)$ if
it is a loop.

These two observations, together with Lemma
\ref{lemma:loopisthmus}, complete the proof of \emph{(i)-(iii)}.
Also, \emph{(iv)} is clear from the definition. $\Box$

\medskip

\begin{defin}
Two matroids $(S_1, \C_1, r_{\C_1})$ and $(S_2, \C_2, r_{\C_2})$
are \emph{isomorphic} if there is a bijection $f:S_1 \rightarrow
S_2$ which induces an isomorphism of simplicial complexes $f:\C_1
\rightarrow \C_2$ such that $r_{\C_1}(c) = r_{\C_2}(f(c))$ for
all $c \in \C_1$.
\end{defin}

A function $f$ on the class $\mathbb{S}$ of semimatroids is called
a \emph{ semimatroid invariant} if $f(\C_1) = f(\C_2)$ for all
$\C_1 \cong \C_2$. An invariant is called a
\emph{Tutte-Grothendieck invariant} (or \emph{T-G invariant}) if
it satisfies the conditions of Proposition \ref{pr:del-cont}. The
following theorem shows that the Tutte polynomial is not only a
T-G invariant; in fact it is the universal T-G invariant on the
class of semimatroids. Any other \emph{generalized T-G invariant},
that is, an invariant satisfying the conditions of Theorem
\ref{th:T-G}, is an evaluation of the Tutte polynomial. An
equivalent result is essentially known for matroids \cite{Br72},
\cite{Ox79}.

\begin{defin}
For a semimatroid $\C = (S, \C, r)$, let $\#\C$ be the number of
elements $x \in S$ such that $\{x\} \in \C$. A semimatroid is
\emph{non-trivial} if $\#\C \neq 0$.
\end{defin}

\begin{theorem}\label{th:T-G}
Let $\mathbb{S}$ be the class of non-trivial semimatroids. Let
$\kk$ be a field and $a, b \in \kk$; and let $R$ be a commutative
ring containing $\kk$. Let $f: \mathbb{S} \rightarrow R$ be a
\emph{generalized T-G invariant}; \emph{i.e.},
\begin{enumerate}
\item[(i)]
if $\C_1 \cong \C_2$ then $f(\C_1) = f(\C_2)$.
\item[(ii)]
If $e \in S$ is neither an isthmus nor a loop in $\C = (S, \C,
r_{\C})$ and $\{e\} \in \C$, then $f(\C) = a f(\C-e) + b f(\C/e).$
\item[(iii)]
If $e$ is an isthmus in $\C$, then $f(\C) = f(I)f(\C-e).$
\item[(iv)]
If $e$ is a loop in $\C$, then $f(\C) = f(L)f(\C/e).$
\item[(v)]
If $e \in S$ and $\{e\} \notin \C$ then $f(S, \C, r_{\C}) = f(S-e,
\C, r_{\C})$.
\end{enumerate}
Then the function $f$ is given by $f(\C) =
a^{\#\C-r_{\C}}\,b^{r_{\C}} \,T_{\C}(f(I)/b, f(L)/a)$ for $\C =
(S, \C, r_{\C})$.

Here $I = (\{i\}, \{\emptyset, \{i\}\}, r)$ denotes the
semimatroid consisting of a single isthmus $i$, and $L (\{l\},
\{\emptyset, \{l\}\}, r)$ denotes the semimatroid consisting of a
single loop $l$.
\end{theorem}

\noindent \emph{Proof.} We can proceed by induction. The only
non-trivial semimatroids which cannot be decomposed using $(ii),
(iii), (iv)$ and $(v)$ are $I$ and $L$, in which case the formula
for $f(\C)$ holds trivially. It simply remains to show that
$a^{\#\C-r_{\C}} \, b^{r_{\C}} T_{\C} (f(I)/b, f(L)/a)$ satisfies
the relations $(ii), (iii), (iv)$ and $(v)$. This is
straightforward from Proposition \ref{pr:del-cont}. $\Box$

\medskip

We conclude this section with some remarks about the relationship
between the Tutte polynomial of a semimatroid $\C$, the Tutte
polynomials of its associated triple $(\Nh, N, M)$, and the Tutte
polynomial of the dual semimatroid $\C^*$.

In the study of the characteristic polynomial $\chi(q)$ of an
affine hyperplane arrangement, the coning construction of
Definition \ref{def:cone} is fundamental, due to the following
result.

\begin{proposition}\emph{(\cite[Proposition 2.51]{Or92})}\label{pr:cone}
For any arrangement $\A$,
\begin{equation*}
\chi_{c\A}(q) = (q-1)\chi_{\A}(q).
\end{equation*}
\end{proposition}

This proposition tells us that, to study characteristic
polynomials of arrangements, we can essentially focus our
attention on central arrangements. Proposition \ref{pr:cone}
generalizes immediately to semimatroids.

As we saw in Theorem \ref{th:pointedmatsemimat}, the analog of
the cone of an arrangement $\A$ is the matroid $\Nh$ of the
semimatroid $\C$. If, in analogy with the definition for
arrangements, we define the \emph{characteristic polynomial of the
semimatroid} $\C$ to be $\chi_{\C}(q) = (-1)^r T_{\C}(1-q,0)$, we
have the following proposition.

\begin{proposition}\label{prop:charsemi}
For any semimatroid $\C$,
\begin{equation*}
\chi_{\Nh}(q) = (q-1)\chi_{\C}(q).
\end{equation*}
\end{proposition}

We might wonder if this result generalizes to the Tutte
polynomial. It turns out that this situation is not so simple. Let
\begin{equation} \label{eq:U}
U_{\C}(x,y) = \sum_{X \notin \, \C}
(x-1)^{r_M-r_M(X)}(y-1)^{|X|-r_M(X)}.
\end{equation}

Then, by looking at the defining sums of $T_M, T_N$ and
$T_{\Nh}$, it is easy to see that $T_M = T_{\C} + U_{\C}$, $T_N =
(x-1)T_{\C} + U_{\C}/(y-1)$, and $T_{\Nh} = xT_{\C} + y/(y-1)
U_{\C}$. (The third of these equations proves Proposition
\ref{prop:charsemi}.) This means that we \emph{can} express the
Tutte polynomial of $\C$ in terms of the Tutte polynomials of
these three matroids $M,N$ and $\Nh$, by solving for $T_{\C}$ in
any two of these three equations. However, $T_{\C}$ does not only
depend on $T_{\Nh}$. This simple dependence takes place for the
characteristic polynomial only because the second term in the
expression of $T_{\Nh}$ vanishes when we substitute $x=1-q$ and
$y=0$.

We conclude that the Tutte polynomial of a semimatroid is closely
related to the Tutte polynomials of its associated triple $(\Nh,
N, M)$. However, the relationship is not simple enough that we
can derive our results on Tutte polynomials of semimatroids as
simple consequences of the analogous results for matroids.

\medskip

Now let us discuss duality and the Tutte polynomial. For matroids
$M$, we know that $T_{M^*}(x,y) = T_M(y,x)$. This is not the case
for a semimatroid $\C$. In fact, it is not difficult to see that
$T_{\C^*}(x,y) = U_{\C}(y,x)/(x-1)$.

It is possible to define a three-variable Tutte-like polynomial of
a semimatroid which is more compatible with duality. In a slightly
different language, this was done by Las Vergnas \cite{La99}, who
defined the concept of the \emph{Tutte polynomial of a quotient
map}. In fact, if the semimatroid $\C$ corresponds to the quotient
map $N \rightarrow M$ under Corollary \ref{th:semimatmatpreim},
then our definition of the Tutte polynomial of $\C$ coincides with
the coefficient of $z$ in Las Vergnas's definition of the Tutte
polynomial of the quotient map $N \rightarrow M$. In particular,
the upcoming Theorem \ref{th:intext} can be derived from his
analogous theorem for quotient maps. His argument uses the
deletion-contraction relation; our approach will give us
additional information about the structure of a semimatroid.

\section{Basis activity.}\label{sec:intext}

We now show that the Tutte polynomial of a semimatroid has
nonnegative coefficients, by giving a combinatorial interpretation
of them. Crapo showed that the coefficients of the Tutte
polynomial of a matroid count the bases with a given internal and
external activity \cite{Cr69}. Our interpretation in the case of
semimatroids is analogous. There are some subtleties involved in
extending this result to semimatroids, so we will need to give
slightly different definitions of internal and external activity.
Our proof will be slightly different from his as well.

In this section we will work with a fixed semimatroid $\C = (S,
\C, r)$. We will denote elements of $S$ by lower case letters,
and subsets of $S$ by upper case letters. As mentioned after
Definition \ref{def:semi}, we will sometimes call the sets in
$\C$ central sets.  Proposition \ref{prop:semi-mat} shows that the
rank function $r$ extends to a matroid rank function on $2^S$,
which we will also call $r$. No confusion arises from this
notation because the semimatroid and matroid rank functions have
the same value where they are both defined.

A \emph{basis} of $\C = (S, \C, r)$ is a set $B \in \C$ such that
$|B| = r(B) = r$. A set $X \in \C$ is \emph{dependent} if $r(X) <
|X|$ and \emph{independent} otherwise. A \emph{circuit} $C$ of
$\C$ is a minimal dependent set in $\C$. Clearly such a set
satisfies $r(C) = |C|-1$. A \emph{cocircuit} $D$ is a minimal
subset of $S$ whose deletion from $\C$ makes the rank of $\C$
decrease; \emph{i.e.}, one such that $r(S-D) < r$, where $r=r(S)$
is the rank of $\C$. Clearly a cocircuit satisfies $r(S-D) = r-1$.

\begin{lemma}\label{l:circuit}
Let $B$ be a basis of $\C$, and let $e \notin B$ be such that $B
\cup e \in \C$. Then $B \cup e$ contains a unique circuit.
\end{lemma}

\noindent \emph{Proof.} Since $B \cup e \in \C$ is dependent, it
contains a circuit. Now assume that it contains two different
circuits $C_1$ and $C_2$. By (R3) we know that
\begin{eqnarray*}
r(C_1 \cap C_2) + r(C_1 \cup C_2) &\leq& r(C_1)+r(C_2) \\
& = & |C_1| -1 + |C_2| -1 \\
& = & |C_1 \cap C_2|-1 + |C_1 \cup C_2| -1.
\end{eqnarray*}
But $r(B \cup e) = |B \cup e| -1$ so, by (R2'), $r(X) \geq |X|-1$
for all $X \subseteq B \cup e$. Therefore $r(C_1 \cap C_2) = |C_1
\cap C_2| -1$ and $r(C_1 \cup C_2) = |C_1 \cup C_2| -1$. Thus
$C_1 \cap C_2$ is a dependent set in $\C$, and it is a proper
subset of the circuit $C_1$. This is a contradiction. $\Box$

\medskip

\begin{lemma}\label{l:cocircuit}
Let $B$ be a basis of $\C$, and let $i \in B$. Then $S - B \cup i$
contains a unique cocircuit.
\end{lemma}

\noindent \emph{Proof.} The deletion of $S - B \cup i$ from $\C$
makes the rank of $\C$ decrease, so this set contains a cocircuit.
Assume that it contains two different cocircuits $B_1$ and $B_2$.
Then
\begin{eqnarray*}
r(S - (B_1 \cap B_2)) & = & r((S - B_1) \cup (S - B_2)) \\
& \leq & r(S-B_1) + r(S-B_2) - r((S - B_1) \cap (S - B_2)) \\
& = & (r - 1) + (r - 1) - r(S-(B_1 \cup B_2)).
\end{eqnarray*}
But $S-(B_1 \cup B_2) \supseteq B-i\, $ and $\,r(B-i) = r-1$, so
$r(S-(B_1 \cup B_2)) \geq r-1$. It follows that $r(S - (B_1 \cap
B_2)) \leq r-1$. Hence the removal of $B_1 \cap B_2$ makes the
rank of the semimatroid decrease, and $B_1 \cap B_2$ is a proper
subset of the cocircuit $B_1$. This is a contradiction. $\Box$

\medskip

From now on, we will fix a linear order on $S$. Now each
$k$-subset of $S$ corresponds to a strictly increasing sequence of
$k$ numbers between $1$ and $|S|$. For each $0 \leq k \leq |S|$,
order the $k$-subsets of $S$ using the lexicographic order on
these sequences.

\begin{defin}
Let $B$ be a basis of $\C$. An element $e \notin B$ is an
\emph{externally active element} for $B$ if $B \cup e \in \C$ and
$e$ is the smallest element\footnote{according to the fixed
linear order} of the unique circuit in $B \cup e$.  Let $E(B)$ be
the set of externally active elements for $B$, and let $e(B) =
|E(B)|$. We call $e(B)$ the \emph{external activity} of $B$.
\end{defin}

\begin{defin}
Let $B$ be a basis of $\C$. An element $i \in B$ is an
\emph{internally active element} in $B$ if $i$ is the smallest
element of the unique cocircuit in $S - B \cup i$. Let $I(B)$ be
the set of internally active elements for $B$, and let $i(B) =
|I(B)|$. We call $i(B)$ the \emph{internal activity} of $B$.
\end{defin}

Now we are in a position to state the main theorem of this
section.

\begin{theorem}\label{th:intext}
For any semimatroid $\C$,
$$
T_{\C}(q,t) = \sum_{B \, \mathrm{basis \,\, of} \, \C} q^{i(B)}
t^{e(B)}.
$$
\end{theorem}

Theorem \ref{th:intext} shows that the coefficients of the Tutte
polynomial are nonnegative integers. The coefficient of $q^i t^e$
is equal to the number of bases of $\C$ with internal activity
$i$ and external activity $e$.

\medskip

We still have some work to do before we can prove Theorem
\ref{th:intext}. The next step will be to give a very useful
characterization of internally and externally active elements.
From now on, when proving results about internally and externally
active elements, we will always use Lemmas \ref{l:ext} and
\ref{l:int} instead of the original definitions.

Given $X \subseteq S$ and an element $e$, let $X_{>e} = \{x \in X
\, | \, x > e\}$. Define $X_{<e}$ analogously.

\begin{lemma}\label{l:ext}
Let $B$ be a basis of $\C$ and let $e \notin B$ be such that $B
\cup e \in \C$. Then $e$ is externally active for $B$ if and only
if $r(B_{>e} \cup e) = r(B_{>e})$.
\end{lemma}

\noindent \emph{Proof.} First assume that $r(B_{>e} \cup e) =
r(B_{>e})$. Then $B_{>e} \cup e \in \C$ is dependent, so it
contains a circuit $C$; $e$ is clearly the smallest element in
this circuit. But $C$ must also be the unique circuit contained
in $B \cup e$. Therefore $e$ is an externally active element for
$B$.

Now assume that $e$ is externally active for $B$. The unique
circuit in $B \cup e$ obviously contains $e$; call it $C \cup e$.
Then $C \subseteq B_{>e}$. By submodularity, we have $r(B_{>e}) +
r(C \cup e) \geq r(B_{>e} \cup e) + r(C)$. But $r(C \cup e) =
r(C)$, so $r(B_{>e}) \geq r(B_{>e} \cup e)$ and the desired
result follows. $\Box$

\medskip

\begin{lemma}\label{l:int}
Let $B$ be a basis and $i \in B$. Then $i$ is internally active in
$B$ if and only if $r(B - i \cup S_{<i}) < r$.\footnote{In fact,
this is true if and only if $r(B - i \cup S_{<i}) = r-1$.}
\end{lemma}

\noindent \emph{Proof.} First assume that $r(B - i \, \cup S_{<i})
< r$. Then the removal of $(S - B)_{> i} \cup i$ makes the rank of
the semimatroid drop, so $(S - B)_{> i} \cup i$ contains a
cocircuit. This cocircuit must contain $i$; call it $D \cup i$,
where $D \subseteq (S - B)_{> i}$. The smallest element of this
cocircuit is $i$, and this cocircuit must also be the unique
cocircuit contained in $S - B \cup i$. Therefore $i$ is an
internally active element of $B$.

Now assume that $i$ is internally active in $B$. Let $S - D \cup
i$ be the unique cocircuit in $S - B \cup i$, where $D \supseteq
B$. Since $i$ is the smallest element in this cocircuit, $D
\supseteq S_{<i}$. Therefore $B \cup S_{<i} \subseteq D$ and,
since $S - D \cup i$ is a cocircuit, $r(B -i \, \cup S_{<i}) <
r(D-i) < r$. $\Box$

\medskip

Now we wish to present a different description of sets in $\C$.
To do it, we need two definitions. For each $X \subseteq S $, let
$dX$ be the lexicographically largest basis of $X$. For each
independent set $X$, which is necessarily in $\C$, let $uX$ be the
lexicographically smallest basis of $\C$ which contains
$X$.\footnote{We will extend the definition of $uX$ to all $X
\subseteq S$ after the proof of Lemma \ref{l:ud}. For simplicity,
we postpone the full definition until then.} Notice that, for any
$X \subseteq S$, $udX$ is a basis of $\C$.

\begin{defin}
Let $\T$ be the set of triples $(B,I,E)$ such that $B$ is a basis
of $\C$, $I \subseteq I(B)$ is a set of internally active
elements for $B$, and $E \subseteq E(B)$ is a set of internally
active elements of $B$.
\end{defin}

We will establish a bijection between $\T$ and $\C$. Define two
maps $\phi_1$ and $\phi_2$ as follows. Given $(B,I,E) \in \T$,
let $\phi_1(B,I,E) = B - I \cup E$. Given $X \in \C$, let
$\phi_2(X) = (udX, udX - dX, X - dX)$. We will show that the maps
$\phi_1$ and $\phi_2$ give the desired bijection: every set $X
\in \C$ can be written uniquely in the form $X = B - I \cup E$
where $B$ is a basis of $\C$, $I \subseteq I(B)$ and $E \subseteq
E(B)$.

\medskip

\noindent \emph{Example.} Recall the arrangement $\A$ introduced
at the beginning of Section \ref{sec:Tuttesemi}. Table
\ref{table:intext} illustrates the bijection between $\T$ and $\C$
in that case. Theorem \ref{th:intext} and Table
\ref{table:intext} imply that $T_{\A}(q,t) = q^3 + 2q^2 +qt + q +
t$, confirming our computation at the beginning of Section
\ref{sec:Tuttesemi}.

\begin{table}[t]
\centering \setlength{\extrarowheight}{2pt}
\begin{tabular}{|c|c|c|l|}
\hline
$B$ & $I(B)$ & $E(B)$ & possible $B-I \cup E$  \\
\hline
$123$ & $123$ & - & $\emptyset, 1, 2, 3, 12, 13, 23, 123$ \\
\hline
$124$ & $12$ & - & $4,14, 24, 124$ \\
\hline
$134$ & $1$ & 2 & $34, 134, 234, 1234$ \\
\hline
$235$ & 23  & $-$ & $5, 25, 35, 235$ \\
\hline
$245$ & 2 & - & $45, 245$ \\
\hline
$345$ & - & 2 & $345, 2345$ \\
\hline
\end{tabular}
\caption{The bijection between $\T$ and $\C$.}
\label{table:intext}
\end{table}

\medskip

\begin{lemma}\label{l:phi1}
The map $\phi_1$ maps $\T$ to $\C$.
\end{lemma}

\noindent \emph{Proof.} Let $(B,I,E) \in \T$. For all $e \in E$,
$B \cup e$ is central and $r(B \cup e) = r(B)$, so $e \in \cl(B)$.
Therefore $E \subseteq \cl(B)$ and $B \cup E \subseteq \cl(B)$.
Since $\cl(B) \in \C$, this implies that $B \cup E \in \C$, and
$B - I \cup E \in \C$ as well. $\Box$

\medskip

\begin{lemma}\label{l:phi2}
The map $\phi_2$ maps $\C$ to $\T$.
\end{lemma}

\noindent \emph{Proof.} Let $X \in \C$. Let $D = dX$ and $U =
udX$, so that $\phi_2(X) = (U, U-D, X-D)$. We need to show three
things.

First, we need $U$ to be a basis for $X$. This is immediate.

Next, we need the elements of $U-D$ to be internally active in
$U$. Let $x \in U-D$. Since $U$ is the smallest basis for $\C$
containing $D$, for any element $x'<x$ not in $U$ we have $r(U-x
\cup x') = r-1 = r(U-x)$ . By submodularity, we can conclude that
$r(U-x \cup S_{<x}) = r-1$, which is exactly what we wanted.

Finally, we need to show that the elements of $X-D$ are externally
active in $U$. Let $x \in X-D$. First notice that $x \notin U$,
because $D \cup x$ is dependent: $r(D \cup x) \leq r(X) = r(D)$.
Also notice that $U \cup x$ is central, applying (CR1) to $D \cup
x$ and $U$. Now observe the following. We know that $D$ is the
largest basis for $X$. Therefore $r(D - x' \cup x) = r(D) - 1$ for
all $x' \in D_{<x}$. By submodularity, it follows that $r(D -
D_{<x} \cup x) = r(D) - |D_{<x}|$. We can rewrite this as
$r(D_{>x} \cup x) = r(D_{>x})$ since $D$ is independent. Since
$D_{>x} \subseteq U_{>x}$, submodularity implies that $r(U_{>x}
\cup x) = r(U_{>x})$. This shows that $x$ is an externally active
element in $U$. $\Box$

\medskip

\begin{proposition}\label{pr:key}
The map $\phi_1$ is a bijection from $\T$ to $\C$, and the map
$\phi_2$ is its inverse.
\end{proposition}

Proposition \ref{pr:key} is the main ingredient of our proof of
Theorem \ref{th:intext}. Before proving it, we need some lemmas.

\begin{lemma}\label{l:rank}
For all $(B,I,E) \in \T$, we have $r(B - I \cup E) = r - |I|$.
\end{lemma}

\noindent \emph{Proof.} We start by showing that $r(B - \, i \,
\cup \, e) = r - 1$ for all $i \in I(B), e \in E(B)$. If $e < i$,
do the following. Since $i$ is internally active, $r(B - \, i \,
\cup S_{<i}) = r-1 = r(B - \, i)$, and therefore $r(B - \, i \,
\cup e) = r-1$. Otherwise, if $i < e$, then $B_{>e} \subseteq B -
i$. Since $e$ is externally active, $r(B_{>e} \cup e) =
r(B_{>e})$. Submodularity then implies that $r(B - \, i \, \cup
e) = r(B-i) = r - 1$.

Now that we know this, submodularity implies that $r(B - i \cup E)
= r-1$ for all $i \in I(B), E \subseteq E(B)$. Applying
submodularity again, we get $r(B - I \cup E) = r - |I|$ for all
$I \subseteq I(B), E \subseteq E(B)$. $\Box$

\medskip

\begin{lemma}\label{l:d}
For all $(B,I,E) \in \T$, we have $d(B - I \cup E) = B - I$.
\end{lemma}

\noindent \emph{Proof.} Lemma \ref{l:rank} tells us that $B-I$ is
a basis for $B - I \cup E$; we need to show that it is the
largest one. Consider an arbitrary $(r - |I|)-$ subset $X$ of $B -
I \cup E$ with $X > B - I$. We will show that $X$ is not a basis
for $B-I \cup E$.

Let $X = (B - I) - (b_1 \cup \cdots \cup b_k) \cup (e_1 \cup
\cdots \cup e_k)$, where the $b_i$'s are in $B-I$ and the $e_i$'s
are in $E$. Since $X > B-I$ we can assume, without loss of
generality, that $b_1 < e_1, \ldots, e_k$.

From Lemma \ref{l:rank} we know that $r(B - I \cup e_i) = r -
|I|$ for all $1 \leq i \leq k$. Also, as we saw in the proof of
Lemma \ref{l:rank}, having $b_1 \in B$, $e_i \in E(B)$ and $b_1 <
e_i$ implies that $r(B - b_1 \cup e_i) = r - 1$. Combining these
two inequalities and using submodularity, we get that $r(B - I -
b_1 \cup e_i) = r - |I| - 1$ for all $1 \leq i \leq k$. Invoking
submodularity once again, we get that $r((B - I) - b_1 \cup (e_1
\cup \cdots \cup e_k)) = r - |I| - 1$. Therefore $r(X) = r((B - I)
- (b_1 \cup \cdots \cup b_k) \cup (e_1 \cup \cdots \cup e_k))
\leq r - |I| - 1 < r(B - I \cup E).$ It follows that $X$ is not a
basis for $B - I \cup E$. $\Box$

\medskip

\begin{lemma} \label{l:ud}
For all $(B,I,E) \in \T$, we have $ud(B - I \cup E) = B$.
\end{lemma}

\noindent \emph{Proof.} In view of Lemma \ref{l:d}, we need to
show that $u(B-I) = B$. Clearly $B$ is a basis of $\C$ containing
$B-I$; now we show that it is the smallest one.

Let $X = B - (b_1 \cup \cdots \cup b_k) \cup (c_1 \cup \cdots \cup
c_k)$ be an $r$-tuple smaller than $B$, where the $b_i$'s are in
$I$ (since $X$ must contain $B-I$) and the $c_i$'s are in $S$. We
will show that $X$ is not a basis for $\C$. Once again we can
assume, without loss of generality, that $c_1 < b_1, \ldots, b_k$.

Since each $b_i$ is internally active, $r(B - b_i \cup S_{<b_i})
= r-1$, and hence $r(B - \, b_i \, \cup c_1) = r-1$.
Submodularity gives $r(B - (b_1 \cup \cdots \cup b_k) \cup c_1) =
r - k$, which in turn gives $r(X) = r(B - (b_1 \cup \cdots \cup
b_k) \cup (c_1 \cup \cdots \cup c_k)) \leq (r-k) + (k-1) < r$.
$\Box$

\medskip

So far we have only defined $uX$ for independent sets $X$ of $\C$.
We can extend the definition to arbitrary subsets $X \subseteq S$
as follows. If $X$ is dependent, then there is no basis of $\C$
containing it. Instead, we consider all the minimal sets of rank
$r$ which contain $X$. Let $uX$ be the lexicographically smallest
of those sets. Then we can say even more.

\begin{lemma}\label{l:du}
For all $(B,I,E) \in \T$, we have $u(B - I \cup E) = B \cup E$ and
$du(B - I \cup E) = B$.
\end{lemma}

We will not need Lemma \ref{l:du} to prove Proposition
\ref{pr:key} and Theorem \ref{th:intext}. We state it for
completeness, but we omit its proof, which is very similar to the
proofs of Lemmas \ref{l:d} and \ref{l:ud}.

\medskip

\noindent \emph{Proof of Proposition \ref{pr:key}.} Checking that
$\phi_1 \circ \phi_2$ is the identity map in $\C$ is immediate,
and Lemmas \ref{l:d} and \ref{l:ud} imply that $\phi_2 \circ
\phi_1$ is the identity map in $\T$. $\Box$

\medskip

\noindent \emph{Proof of Theorem \ref{th:intext}.} Using the
bijection of Proposition \ref{pr:key}, the sets in $\C$ are
precisely the sets of the form $B-I \cup E$, where $B$ is a basis,
$I \subseteq I(B)$ and $E \subseteq E(B)$. Also, from Lemma
\ref{l:rank}, $r(B-I \cup E) = r - |I|$. Therefore we have
\begin{eqnarray*}
T(q,t) &=&  \sum_{X \in \, \C}
(q-1)^{r-r(X)}(t-1)^{|X|- r(X)}  \\
&=& \sum_{B \, \mathrm{basis}}\sum_{I \subseteq I(B)}\sum_{E
\subseteq E(B)} (q-1)^{r- r(B - I\, \cup E)} (t-1)^{|B - I\, \cup
E| - r(B - I\,
\cup E)} \\
&=& \sum_{B\, \mathrm{basis}}\sum_{I \subseteq I(B)}\sum_{E
\subseteq E(B)}
(q-1)^{|I|}(t-1)^{|E|} \\
&=& \sum_{B\, \mathrm{basis}}(1+(q-1))^{|I(B)|}(1+(t-1))^{|E(B)|} \\
&=& \sum_{B\, \mathrm{basis}} q^{i(B)}t^{e(B)}.
\end{eqnarray*}
as desired. $\Box$

\medskip

\begin{figure}[!ht]
\centering
\includegraphics[width=4.5in]{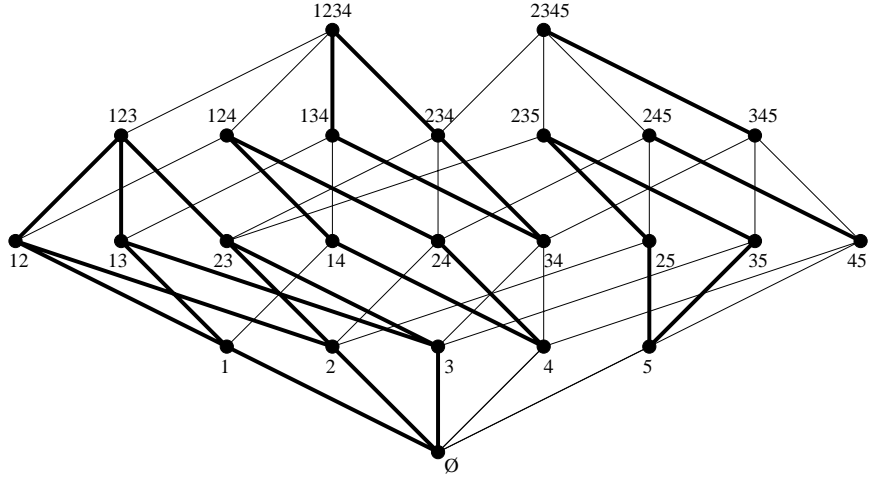}
\caption{The decomposition of $\C$ into intervals.}
\label{fig:Boolean}
\end{figure}

Regard the simplicial complex $\C$ as a poset, ordering its faces
by inclusion. There is a nice way to understand Theorem
\ref{th:intext} in terms of this poset. Proposition \ref{pr:key}
gives us a way of classifying the faces of $\C$ according to the
basis of $\C$ that they correspond to under the map $ud$ (or
$du$). This classification decomposes the poset into disjoint
intervals, where each interval is a Boolean algebra of the form
$[B - I(B), B \cup E(B)]$ for a basis $B$. This is illustrated in
Figure \ref{fig:Boolean} for the arrangement $\A$ considered at
the beginning of Section \ref{sec:Tuttesemi}; recall Table
\ref{table:intext}. If we look at the interval corresponding to
basis $B$, and add the contributions of its elements to the
right-hand side of (\ref{eq:Tutte}), we simply get the monomial
$q^{i(B)}t^{e(B)}$.

\section{Acknowledgments.}

The present work is Chapter 3 of the author's Ph.D. thesis
\cite{Ar02}. I am extremely grateful to Seth Chaiken, Vic Reiner,
Gian-Carlo Rota, Richard Stanley, and Tom Zaslavsky; the content
and exposition of this paper benefitted greatly from instructive
conversations with them.

\end{document}